\documentclass[12pt,a4paper]{article}
\usepackage{amssymb,amsmath}
\usepackage{amsfonts}
\usepackage{amsthm}
\usepackage{a4}
\usepackage{tikz}
\usepackage{enumerate}%
\usepackage{verbatim}
\usepackage{mathalfa}
\usepackage{caption}
\usepackage{float}
\newcommand{\para}{\par\vspace{.25cm}}

\newtheorem*{theorem*}{Theorem}
\newtheorem{theorem}{Theorem}[section]

\newtheorem{cor}[theorem]{Corollary}
\newtheorem{remark}[theorem]{Remark}
\newtheorem{notation}[theorem]{Notation}

\newcommand{\Q}{\mathbb{Q}}

\newcommand{\C}{\operatorname{Cen}}

\newcommand{\G}{\operatorname{Gal}}

\begin{document}
		\title{\bf Matrix units in the simple components of rational group algebras}
		\author{ Gurmeet K. Bakshi{\footnote {Research supported by DST-FIST grant no. SR/FST/MS-II/2019/43 is
		gratefully acknowledged.}}  ~and Jyoti Garg{\footnote {Research supported by Council of Scientific and Industrial Research (CSIR), Govt. Of India under the reference no. 09/135(0886)/2019-EMR-I is gratefully acknowledged.} \footnote{Corresponding author}} \\ {\em \small Department of 
				Mathematics,}\\
			{\em \small Panjab University, Chandigarh 160014, India}\\{\em
				\small email: gkbakshi@pu.ac.in and jyotigarg0811@gmail.com} }
		\date{}
		{\maketitle}
\begin{abstract}
For the rational group algebra $\Q G$ of a finite group $G$, we provide an effective method to compute a complete set of matrix units and, in particular, primitive orthogonal idempotents in a simple component of $\Q G$, which is realized by a generalized strongly monomial character and has a prime Schur index. We also provide some classes of groups $G$ where this method can be successfully applied. The application of the method developed is also illustrated with detailed computations.
\end{abstract}	
\noindent\textbf{Keywords}: Shoda pairs, generalized strong Shoda pairs, primitive idempotents, matrix units, Schur index, Frobenius groups.
\para \noindent {\bf MSC2000:} 17C27, 20C05, 16K20, 16S35, 16U40.
\section{Introduction}
Let $G$ be a finite group and $\Q G$ the rational group algebra of $G$. The knowledge of a complete set of primitive orthogonal idempotents and more generally matrix units of $\Q G$ is useful tool for various problems concerning
\begin{itemize}
\item a full description of the Wedderburn decomposition of $\Q G$;
\item automorphism group of the rational group algebras;
\item rational representations of $G$;
\item unit group of $\mathbb{Z}$-orders in $\Q G$ and, in particular, the unit group of integral group ring $\mathbb{Z}G.$
\end{itemize}\noindent Let us begin by describing the notion of primitive orthogonal idempotents and matrix units of $\Q G$. Let $\operatorname{Irr}(G)$ be the set of all complex irreducible characters of $G$. For $\chi\in\operatorname{Irr}(G)$, $${e_{\mathbb{Q}}(\chi) = \frac{\chi(1)}{|G|}\sum_{\sigma \in \operatorname{Gal}(\mathbb{Q}(\chi)/ \mathbb{Q})} {\sum _{g \in G}} \sigma(\chi(g))g^{-1}}$$ is a  primitive central idempotent of $\Q G$. Here $|G|$ is the order of $G$, $\Q(\chi)$ is the field obtained by adjoining to $\Q$  the character values $\chi(g),~g \in G$, and $\G(\Q(\chi)/\Q)$ is the Galois group of extension $\Q(\chi)/\Q.$ This implies that each $\chi\in\operatorname{Irr}(G)$ give arise to the simple component $\Q Ge_{\Q}(\chi)$ of $\Q G$. Moreover, every simple component of $\Q G$ arises this way. For $\chi\in\operatorname{Irr}(G)$, the simple component $\Q Ge_{\Q}(\chi)$ is well known to be isomorphic to a matrix algebra $M_n(D)$ over a division algebra $D$ containing $\Q$ in its center. Furthermore, it has a unique set of elements $\{E_{ij}~|~1\leq i,j\leq n\},$ called a complete set of matrix units of $\Q Ge_{\Q}(\chi)$ such that $e_{\Q}(\chi)=\sum_{i=1}^{n}E_{ii}$ and $E_{ij}E_{kl}=\delta_{jk}E_{il}$ for all $i,j,k,l$, where $\delta_{jk}$ is Kroneckar delta function. The particular matrix units $\{E_{ii}~|~1\leq i\leq n\}$ is called a complete set of primitive orthogonal idempotents of $\Q Ge_{\Q}(\chi)$. The union of complete sets of matrix units (respectively primitive orthogonal idempotents) in all the simple components of $\Q Ge_{\Q}(\chi)$, as $\chi$ runs over $\operatorname{Irr}(G)$, is called a complete set of matrix units (respectively primitive orthogonal idempotents) of $\Q G$.\para \noindent For a nilpotent group $G$, Jespers, Olteanu and del R\'io \cite{JOdR12} gave a complete set of primitive orthogonal idempotents and matrix units of $\Q G$ and hence concluded the earlier result of Roquette \cite{Roq}, which says that any simple components of $\Q G,~G$ nilpotent, must have the Schur index at most $2$. In \cite{JOdRVG13}, Jespers, Olteanu, del R\'io and Van Gelder gave a complete set of primitive orthogonal idempotents and matrix units of $\Q Ge_{\Q}(\chi)$ when $\chi$ is a strongly monomial irreducible character such that the Schur index of $\Q Ge_{\Q}(\chi)$ is $1$. In a recent joint work with Gabriela Olteanu \cite{BGO}, we extended this result from a strongly monomial character to a generalized strongly monomial character $\chi$, while keeping the constraint that the Schur index of $\Q Ge_{\Q}(\chi)$ is $1$. For the definition of strongly monomial character and generalized strongly monomial character, we refer to section $2$ (Preliminaries). It may be mentioned that the generalized strongly monomial characters form a vast class of monomial characters (see \cite{BK19,BK22,BGK22}).\para\noindent In section $3$, we provide an effective method to compute a complete set of matrix units and, in particular, primitive orthogonal idempotents of $\Q Ge_{\Q}(\chi)$, when $\chi$ is a generalized strongly monomial character such that the Schur index of $\Q Ge_{\Q}(\chi)$ is a prime. As mentioned earlier, $\Q Ge_{\Q}(\chi)$ is isomorphic to $M_n(D)$ for some division algebra $D$. In order to write the matrix units and, in particular, primitive orthogonal idempotents of $\Q Ge_{\Q}(\chi)$, the key hurdle is to understand a split $F$-subalgebra of $\Q Ge_{\Q}(\chi)$ isomorphic to $M_n(F)$, where $F$ is the center of $\Q Ge_{\Q}(\chi)$, whose centralizer in $\Q Ge_{\Q}(\chi)$ is a division algebra. This is precisely what have been achieved in Theorems \ref{l1}, \ref{3thm1} and \ref{3thm3}.\para\noindent In section $4$, we show that the method can be applied to determine a complete set of matrix units and, in particular, primitive orthogonal idempotents of $\Q G$ for the following groups $G$:
\begin{itemize}
\item A class of Frobenius groups and, in particular, nilpotent-by-nilpotent groups with fixed-point-free action.
\item Camina groups.
\item Groups of order $p_1p_2\cdots p_k$ and $p_1^2p_2\cdots p_k$, where $p_{i}$$'$s are distinct primes.
\item Groups of order $p_{1}^{2}p_{2}^2$, where $p_1$ and $p_2$ are primes.
\end{itemize} In an earlier work \cite{BG}, we proved that when $G$ is SmallGroup(1000,86) in GAP library, all the simple components of $\Q G$ have the Schur indices $1$ except one component, where it is $2$. In section $4$, we have also illustrated our method to compute a complete set of matrix units and, in particular, primitive orthogonal idempotents in the simple component of $\Q G$ with the Schur index $2$. 
\section{Preliminaries}
We begin with recalling some basic notions on crossed product algebras. For details on crossed product algebras, see (Section 2.6 of \cite{JdR16}, Section 8.5 of \cite{MS02}).\para\noindent Let $R$ be a ring with unity and $G$ a finite group. Let $\operatorname{Aut}(R)$ and $\mathcal{U}(R)$ denote the group of automorphisms of $R$ and the group of units of $R$ respectively. Let $\sigma: G \rightarrow \operatorname{Aut}(R)$ and $\tau : G \times G \rightarrow \mathcal{U}(R)$ be two maps which satisfy the following relations:
$$\tau_{gh,x}\sigma_{x}(\tau_{g,h})= \tau_{g,hx} \tau_{h,x}~~{\rm and} ~~\tau_{g,h}\sigma_{h}(\sigma_{g}(r))= \sigma_{gh}(r) \tau_{g,h},$$for all $g,h,x \in G$ and $r \in R$. Here $\sigma_{g}$ is the image of $g$ under the map $\sigma$ and $\tau_{g,h}$ is the image of $(g,h)$ under the map $\tau.$ \para\noindent Let $R*_{\tau}^{\sigma} G$ denote the set of finite formal sums $\left\{ \sum z_{g}a_{g} ~|~ a_{g} \in R, g \in G \right\}$, where $z_{g}$ is a symbol corresponding to $g \in G.$ Equality and addition in $R*_{\tau}^{\sigma} G$ are defined componentwise. For $g,h \in G$ and $r \in R,$ by setting
$$z_g z_h = z_{gh} \tau_{g,h},~rz_g = z_g \sigma_g(r)$$ and extending this rule distributively, $R*_{\tau}^{\sigma}G$ becomes an associative ring, called the {\it crossed product} of $G$ over $R$ with factor set $\tau$ and action $\sigma$.\para\noindent Note that $R*_{\tau}^{\sigma}G$ is a free $R$-module with $\{z_{g}~|~g \in G\}$ an $R$-basis, which we call a basis of units of $R*_{\tau}^{\sigma} G$ as an $R$-module.\para\noindent A \textit{classical crossed product} is a crossed product $L*_{\tau}^{\sigma} G$, where $L$ is a finite Galois extension of the center  $F$ of $L*_{\tau}^{\sigma}G$, $G=\G(L/F)$ is the Galois group of $L/F$ and $\sigma$ is the natural action of $G$ over $L$. Such a classical crossed product $L*_{\tau}^{\sigma} G$ is commonly denoted by $(L/F,\tau).$\para\noindent By an {\it abelian crossed product algebra}, we mean a classical crossed product $(L/F,\tau)$, where $\operatorname{Gal}(L/F)$ is abelian.\para\noindent A \textit{cyclotomic algebra} over $F$ is a classical crossed product $(L/F,\tau),$ where $L=F(\zeta)$ is a cyclotomic extension over $F$. Here $\zeta$ is a root of unity. \para\noindent 
If $A$ is a central simple $F$-algebra and $E$ is a subfield of $A$ with the property that $\operatorname{dim}_{F}(A)$, the dimension of $A$ over $F$, is equal to $[E:F]^2$, then we call $E$ as a {\it strictly maximal subfield} of $A$. Here, $[E:F]$ denote the degree of the extension $E$ over $F$. Observe that if $A=(E/F,\tau)$ is a classical crossed product, then $E$ is a strictly maximal subfield of $A$. Conversely, we will show in the next paragraph that if a central simple algebra $A$ has a strictly maximal subfield $E$ such that $E/F$ is a Galois extension, then $A$ is isomorphic to $(E/F,\tau)$ for some factor set $\tau.$ 
\para\noindent Let $G=\{\sigma_1,\sigma_2,\cdots,\sigma_n\}$ be the Galois group of $E/F$. For each $\sigma_i \in G,~1\leq i\leq n$, by (\cite{JdR16}, Theorem 2.1.9), there exists an invertible element $u_i\in A$ such that $\sigma_i(x)=u_{i}^{-1}xu_{i}$ for all $x\in E.$ Now, we prove that $u_{i}$$'$s are linearly independent over $E$. Suppose $\sum_{i=1}^{n}u_{i}a_{i}=0$ with $a_i\in E$. Amongst all the expressions of $0$ as $E$ linear combination of $u_1,u_2,\cdots,u_n$, by reordering $u_i$$'$s, suppose $m$ is the least integer such that $\sum_{i=1}^{m}u_{i}a_{i}=0$ with all $a_{i}\neq 0$, $1\leq i\leq m$. For any $x \in E,$ we have that $\sum_{i=1}^{m}u_{i}a_{i}x=0$, which gives $\sum_{i=1}^{m}\sigma_{i}^{-1}(x)u_{i}a_{i}=0$. Consequently, $\sum_{i=1}^{m}\sigma_{i}^{-1}(x)u_{i}a_{i}-\sigma_{m}^{-1}(x)\sum_{i=1}^{m}u_{i}a_{i}=0,$ i.e., $\sum_{i=1}^{m-1}(\sigma_{i}^{-1}(x)-\sigma_{m}^{-1}(x))u_{i}a_{i}=0$. Since $x\in E$ is arbitrary, it is a contradiction to $m$ being the least. This proves that $u_i$$'$s are linearly independent over $E$. Since $\operatorname{dim}_{E}(A)=n,$ we obtain that $A=\{\sum_{i=1}^{n}u_{i}e_{i}~|~ e_{i}\in E\}.$ If $\tau:G\times G\rightarrow \mathcal{U}(E)$ is given by $(u_{i},u_{j})=u_{ij}^{-1}u_{i}u_{j}$, then the associativity of multiplication in $A$ implies that $\tau$ is a factor set and $A$ is isomorphic to the classical crossed product $(E/F,\tau)$.\para\noindent 
A classical theorem due to Brauer and Witt (\cite{JdR16}, Theorem 3.7.1), which connects the simple components of $\Q G$ with the cyclotomic algebras, states that every simple component of the rational group algebra $\mathbb{Q}G$ of a finite group $G$ is Brauer equivalent to a cyclotomic algebra containing $\mathbb{Q}$ in its center.\para\noindent In \cite{BG}, Brauer and Witt's theorem has been improved for the simple component $\Q Ge_{\Q}(\chi)$ of $\Q G$ for certain monomial characters $\chi$. More precisely, it is proved there that for a generalized strongly monomial character $\chi$ (which we shall recall soon), $\Q Ge_{\Q}(\chi)$ is indeed isomorphic to a matrix algebra over a cyclotomic algebra. This result is of key importance in our findings of a complete set of matrix units and, in particular, primitive orthogonal idempotents in this paper. In order to understand its statement, we are now going to recall the notions of Shoda pairs, generalized strong Shoda pairs and generalized strongly monomial characters (groups) introduced in \cite{BK19} and \cite{OdRS04}.\para\noindent  By a Shoda pair $(H,K)$ of $G$, we mean a pair of subgroups of $G$ satisfying the following:
\begin{itemize}
	\item [(i)] $K \unlhd H$,  $H/K$ is cyclic;
	\item [(ii)] if $g \in G$ and $[H, g] \cap H \subseteq K,$ then $ g \in H.$ Here $[H,g]=\langle g^{-1}h^{-1}gh\, |\,  h \in H \rangle.$\end{itemize}
\para \noindent If $(H,K)$ is a Shoda pair of $G$ and $\lambda$ a linear character of $H$ with kernel $K$, then by (\cite{CR62}, Corollary 45.4), we have that $\lambda^G$ is irreducible. We call this irreducible character $\lambda^G$ as an {\it irreducible character of $G$ arising from the Shoda pair $(H,K)$} and the corresponding simple component $\Q Ge_{\Q}(\lambda^G)$ as {\it the simple component of $\Q G$ realized by the irreducible character $\lambda^G$.} Conversely, if $\chi$ is a monomial character of $G$, i.e., $\chi=\lambda^G$, where $\lambda$ is a linear character of a subgroup $H$ of $G$, then $(H,K)$ with $K=\operatorname{ker}(\lambda)$ is a Shoda pair of $G$ (\cite{OdRS04}, Definition 1.4). Hence, the knowledge of monomial irreducible characters of $G$ translate to knowing the Shoda pairs of $G$.
\para \noindent If $(H,K)$ is a Shoda pair of $G$ and $\lambda$ is a linear character of $H$ with kernel $K$, then the pair $(H,K)$ is a called generalized strong Shoda pair of $G$ if there is a chain $H = H_0 \leq H_1 \leq \cdots \leq H_n = G$ (called a strong inductive chain from $H$ to $G$ of length $n$) of subgroups of $G$ such that the following conditions hold for all $0 \leq i \leq n-1$:
\begin{enumerate}[(i)]
\item $H_i \unlhd \operatorname{Cen}_{H_{i+1}}(e_{\Q}(\lambda^{H_i}))$;
\item the distinct $H_{i+1}$-conjugates of $e_{\Q}(\lambda^{H_i})$ are mutually orthogonal.
\end{enumerate}
\noindent Furthermore, if a generalized strong Shoda pair $(H,K)$ has a strong inductive chain of length $1$, then it is called a {\it strong Shoda pair} of $G$ defined by Olivieri, del R\'io and Sim\'on (see \cite{OdRS04}, Definition 3.1).\para\noindent
A finite group $G$  is called {\it generalized strongly monomial group} (respectively {\it strongly monomial group}) if every irreducible character of $G$ arises from a generalized strong Shoda pair (respectively strong Shoda pair) of $G$. The class of generalized strongly monomial groups is a vast class of monomial groups and so far there is no known example of a monomial group which is not generalized strongly monomial group (see \cite{BK19,BK22,BGK22}).\para\noindent In order to state the result from \cite{BG}, which is of key importance, we also need to recall the following notation adopted in that paper and which we shall also continue to use throughout the paper.
\begin{notation}\label{6not1} For a generalized strong Shoda pair $(H,K)$ of $G$, $\lambda$ a linear character of $H$ with kernel $K$ and $$H=H_{0}\leq H_{1}\leq \cdots \leq H_{n}=G$$ a strong inductive chain from $H$ to $G$:
\begin{equation*}
\begin{array}{lll} \noindent
\mathcal{A} &:=&\Q Ge_{\Q}(\lambda^G).\\
C_{i}&:=& \operatorname{Cen}_{H_{i+1}}(e_{\Q}(\lambda^{H_{i}})).\\
k_{i}&:=& [H_{i+1}:C_{i}].\\
k &:=& k_{0}k_{1} \cdots k_{n-1}.\\
\mathtt{k}&:=& \prod_{i=0}^{n-1}|C_{i}/H_{i}|,~\mbox{which is equal to }~\frac{[G:H]}{k}.\\
\mathbb{F} &:=& \mathcal{Z}(\mathcal{A}), \mbox{the center of~} \mathcal{A}.\\
\mbox{E} &:=& \Q He_{\Q}(\lambda).\\
\mathbb{E} &:=& k \times k~ \mbox{scalar matrices}~ \operatorname{diag}(\alpha,\alpha,\cdots,\alpha)_{k}, \alpha \in \mbox{E}.\\
\mbox{F} &:=&\{\alpha \in \mbox{E} ~|~ \operatorname{diag}(\alpha,\alpha,\cdots,\alpha)_{k} \in \mathbb{F}\}.\\
\mathcal{B} &:=& M_{{k}}(\mbox{F}).\\
\operatorname{Cen}_{\mathcal{A}}(\mathcal{B})&:=& \mbox{the centralizer of~} \mathcal{B} \mbox{~in~} \mathcal{A}.\\
 \mathcal{G} &:=& \operatorname{Gal}(\mathbb{E}/\mathbb{F}).
\end{array}
\end{equation*}
\end{notation}\noindent
For a Shoda pair $(H,K)$ of $G$, let $$\varepsilon(H,K)=\left\{\begin{array}{ll}\widehat{K}, & \hbox{$H=K$;} \\\prod(\widehat{K}-\widehat{L}), & \hbox{otherwise,}\end{array}\right.$$ where  $\widehat{H}=\frac{1}{|H|}\displaystyle\sum_{h \in H}h$ and $L$ runs over the normal subgroups of $H$ minimal with respect to the property of including $K$ properly. In (\cite{OdRS04}, Proposition 1.1), it is proved that $\varepsilon(H,K)$ is an idempotent of the group algebra $\mathbb{Q}G$.\para\noindent
We are now ready to state the result from \cite{BG}.
\begin{theorem}{\label{thm2}}{\rm (\cite{BG}, Section 3)}
	Let $(H,K)$ be a generalized strong Shoda pair of $G$, $\lambda$ a linear character of $H$ with kernel $K$ and $H=H_{0}\leq H_{1} \leq \cdots \leq H_{n}=G$ a strong inductive chain from $H$ to $G$. Then 
	\begin{enumerate}[(i)]
	\item $\Q Ge_{\Q}(\lambda^G) \cong \mathcal{B} \otimes_{\mathbb{F}} \operatorname{Cen}_{\mathcal{A}}(\mathcal{B}) \cong M_{k}(\operatorname{Cen}_{\mathcal{A}}(\mathcal{B})).$
	\item $\mathbb{E}$ is a Galois extension of $\mathbb{F}$, $\operatorname{dim}_{\mathbb{F}}(\mathbb{E})=\mathtt{k}=\prod_{i=0}^{n-1}|C_{i}/H_{i}|$, $\operatorname{Cen}_{\mathcal{A}}(\mathcal{B})$ contains $\mathbb{E}$ as a maximal subfield and $\operatorname{dim}_{\mathbb{F}}(\operatorname{Cen}_{\mathcal{A}}(\mathcal{B}))=(\operatorname{dim}_{\mathbb{F}}(\mathbb{E}))^2=\mathtt{k}^2.$
	\item there exist units $\{z_{\sigma}~|~\sigma \in \mathcal{G}\}$ in $\operatorname{Cen}_{\mathcal{A}}(\mathcal{B})$, which form a basis of $\operatorname{Cen}_{\mathcal{A}}(\mathcal{B})$ as a vector space over $\mathbb{E}$ and $\operatorname{Cen}_{\mathcal{A}}(\mathcal{B})\cong (\mathbb{E}/\mathbb{F},\tau),$ where $\tau: \mathcal{G}\times \mathcal{G} \rightarrow \mathcal{U}(\mathbb{E})$ is given by
	\begin{equation}{\label{n}}
	\tau(\sigma,\sigma')=z_{\sigma\sigma'}^{-1}z_{\sigma}z_{\sigma'}.
	\end{equation} 
	\item $\mathcal{G}$ is of the type $(((C_{0}/H_{0}$-by-$C_{1}/H_{1})$-by-$C_{2}/H_{2})\cdots)$-by-$C_{n-1}/H_{n-1}.$ In particular, the order of $\mathcal{G}$ is $\mathtt{k}$. 
	\end{enumerate}
\end{theorem}
\noindent
It may be recalled that a group $G$  is said to be of type $G_{1}$-by-$G_{2}$ if $G$ has a normal subgroup isomorphic to $G_{1}$ and the quotient by $G_{1}$ is isomorphic to $G_{2}$. \begin{quote} In view of Theorem \ref{thm2}, if $\mathcal{E}_{1}$ and $\mathcal{E}_{2}$ denote a complete set of matrix units (primitive orthogonal idempotents) of the $\mathbb{F}$-subalgebras of $\mathcal{A}$ which maps to $\mathcal{B}$ and $\operatorname{Cen}_{\mathcal{A}}(\mathcal{B})$ respectively, under the above isomorphism, then $$\mathcal{E}_{1}\mathcal{E}_{2}=\{e_{1}e_{2}~|~e_{1}\in \mathcal{E}_{1},e_{2}\in\mathcal{E}_{2}\}$$ is a complete set of matrix units (primitive orthogonal idempotents) of $\mathcal{A}$.\end{quote}\para\noindent For the convenience of language, we will denote the central simple $\mathbb{F}$-subalgebra of $\Q Ge_{\Q}(\lambda^G)$ which maps to $\mathcal{B}$ (respectively $\operatorname{Cen}_{\mathcal{A}}(\mathcal{B})$) again by $\mathcal{B}$ (respectively $\operatorname{Cen}_{\mathcal{A}}(\mathcal{B})$).\para \noindent A complete set of primitive orthogonal idempotents of $\mathcal{B}$ has already been computed in Lemma 3.1 of \cite{BGO} (irrespective of the Schur index of $\mathcal{A}$) (see Theorem \ref{l1} stated below). Recall that the {\it Schur index} of a central simple algebra $A$ is defined as the square root of dimension of the division algebra $D$ over its center $F$, once we write $A$ as isomorphic to some $M_n(D)$. 
\begin{theorem}{\label{l1}}{\rm  (\cite{BGO}, Lemma 3.1)}
  Let $(H,K)$ be a generalized strong Shoda pair of $G$, {$\lambda$ a linear character of $H$ with kernel $K$ and} $H=H_{0}\leq H_{1} \leq \cdots \leq H_{n}=G$ a strong inductive chain from $H$ to $G$. For $0 \leq i \leq n-1$, let $T_{i}=\{t_{(i,1)},t_{(i,2)}, \ldots, t_{(i,k_{i})} \}$ be a left transversal of $C_{i}$ in $H_{i+1}$ {with $t_{(i,1)}=1$}, where $C_i=\operatorname{Cen}_{H_{i+1}}(e_{\Q}(\lambda^{H_{i}}))$ and $k_i=[H_{i+1}:C_{i}]$. Let $T= T_{0}T_{1}\cdots T_{n-1}$\linebreak $=\{t_{(0,s_{0})}t_{(1,s_{1})}\cdots t_{(n-1,s_{n-1})}~|~t_{(i,s_{i})} \in T_{i},~0 \leq i \leq n-1,~1 \leq s_{i} \leq k_{i} \}$. Then $$\{ t^{-1} \,\varepsilon(H,K)\,t~|~ t \in T\}$$ is a complete set of primitive orthogonal idempotents of $\mathcal{B}$. 
\end{theorem} 
\begin{remark}\label{rr}
By following the proof of Lemma 3.1 of \cite{BGO}, one can also show that $\{t^{-1} \,\varepsilon(H,K)\,t'~|~ t,t' \in T\}$ is a complete set of matrix units of $\mathcal{B}$.
\end{remark}
\noindent  In (\cite{BGO}, Lemma 3.2 and Theorem 3.5), a complete set of primitive orthogonal idempotents of $\C_{\mathcal{A}}(\mathcal{B})$ and more generally, matrix units of $\mathcal{A}$ is computed under the constraint when $\Q Ge_{\Q}(\lambda^G)$ is a split $\mathbb{F}$-algebra, i.e., when the Schur index of $\Q Ge_{\Q}(\lambda^G)$ is $1$.
\para\noindent The objective of the next section is to move from the Schur index $1$ to a prime Schur index.
\section{Matrix units and primitive idempotents}
Let $G$ be a finite group, $(H,K)$ a generalized strong Shoda pair of $G$ and $\lambda$ a linear character of $H$ with kernel $K$. Let $\mathcal{A}=\Q Ge_{\Q}(\lambda^G)$ be the simple component of $\Q G$ realized by the irreducible character $\lambda^G$. In this section, we will provide a method to compute a complete set of matrix units and primitive orthogonal idempotents of $\mathcal{A}$ in the particular case when this simple component has a prime Schur index.\para\noindent To achieve our objective, the following strategy shall be followed:
\begin{enumerate}[(i)]
\item In Theorem \ref{3thm1}, we will express $\operatorname{Cen}_{\mathcal{A}}(\mathcal{B})$ as a tensor product of two abelian crossed product $\mathbb{F}$-subalgebras, one of which is a split $\mathbb{F}$-subalgebra and the other one has the Schur index $p$ (a prime) and is of dimension (over $\mathbb{F}$) a power of $p$. Hence, a complete set of primitive orthogonal idempotents and matrix units of $\operatorname{Cen}_{\mathcal{A}}(\mathcal{B})$ can be obtained from a complete set of primitive orthogonal idempotents and matrix units of these two subalgebras.
\item In Theorems \ref{3thm2} and \ref{3thm3}, we will explicitly provide a complete set of primitive orthogonal idempotents and matrix units of the following type of crossed product algebras:
\para\noindent$(a)$ an abelian crossed product split $\mathbb{F}$-algebra;\para\noindent $(b)$ an abelian crossed product algebra of the Schur index $p$ and dimension (over $\mathbb{F}$) a power of $p$.   
\end{enumerate}
\begin{theorem}\label{3thm1}
Let $(H,K)$ be a generalized strong Shoda pair of $G$ and $\lambda$ a linear character of $H$ with kernel $K$. Let $\mathcal{A}=\Q Ge_{\Q}(\lambda^G)$ and assume that the Schur index of $\mathcal{A}$ is $p$ (a prime). Let $\mathbb{E},\mathcal{B},\operatorname{Cen}_{\mathcal{A}}(\mathcal{B})=(\mathbb{E}/\mathbb{F},\tau), z_{\sigma_i}$$'$s be as stated above. Let $\mathcal{G}=\operatorname{Gal}(\mathbb{E}/\mathbb{F})$, $\mathcal{G}_{p}$ the Sylow $p$-subgroup of $\mathcal{G}$, $\mathcal{G'}$ the subgroup of $\mathcal{G}$ of order $\frac{|\mathcal{G}|}{|\mathcal{G}_{p}|}$, $\mathbb{E}^{\mathcal{G'}}$ the fixed subfield of $\mathbb{E}$ by the action of $\mathcal{G'}$. The following statements hold:
\begin{enumerate}[(i)]
\item $\tau_{|_{\mathcal{G'}}}$, the restriction of $\tau$ on $\mathcal{G'}\times \mathcal{G'}$, is equivalent to the trivial factor set. The $\mathbb{F}$-subalgebra of $\operatorname{Cen}_{\mathcal{A}}(\mathcal{B})$ generated by $\mathbb{E}^{\mathcal{G}_{p}}$ and $\mathcal{G'}$ is split and is equal to $\mathbb{E}^{\mathcal{G}_{p}}*_{\tau_{|_{\mathcal{G'}}}}\mathcal{G'}$.
\item the centralizer of $\mathbb{E}^{\mathcal{G}_{p}}*_{\tau_{|_{\mathcal{G'}}}}\mathcal{G'}$ in $\operatorname{Cen}_{\mathcal{A}}(\mathcal{B})$ is $\mathbb{E}^{\mathcal{G'}}*_{\tau'}\mathcal{G}_{p}$ for some factor set $\tau'.$ Furthermore, the Schur index of $\mathbb{E}^{\mathcal{G'}}*_{\tau'}\mathcal{G}_{p}$ is $p$ and its dimension over $\mathbb{F}$ is a power of $p$.
\item $\operatorname{Cen}_{\mathcal{A}}(\mathcal{B})\cong \mathbb{E}^{\mathcal{G}_{p}}*_{\tau_{|_{\mathcal{G'}}}}\mathcal{G'}\otimes_{\mathbb{F}}\mathbb{E}^{\mathcal{G'}}*_{\tau'}\mathcal{G}_{p}$. Also, if $\mathcal{E}$ and $\mathcal{E'}$ are complete sets of matrix units (primitive orthogonal idempotents) of $\mathbb{E}^{\mathcal{G}_{p}}*_{\tau_{|_{\mathcal{G'}}}}\mathcal{G'}$ and $\mathbb{E}^{\mathcal{G'}}*_{\tau'}\mathcal{G}_{p}$ respectively, then $$\mathcal{EE'}=\{ee'~|~e\in\mathcal{E},e'\in\mathcal{E'}\}$$ is a complete set of matrix units (primitive orthogonal idempotents) of $\operatorname{Cen}_{\mathcal{A}}(\mathcal{B})$. 
\end{enumerate} 
\end{theorem}
\noindent\textbf{Proof.} $(i)$ As the Schur index of $\operatorname{Cen}_{\mathcal{A}}(\mathcal{B})$ divides $|\mathcal{G}|,$ we have that $p||\mathcal{G}|.$ Let $|\mathcal{G}|=p_1^{n_1}p_2^{n_2}\cdots p_r^{n_r}$, where $p_i$$'$s are distinct primes, $p_1=p$ and $n_{i}\geq 1.$ If $\mathcal{G}_{p}$ is the Sylow $p$-subgroup of $\mathcal{G}$ and $\mathcal{G}'$ is the subgroup of $\mathcal{G}$ of order $\frac{|\mathcal{G}|}{|\mathcal{G}_{p}|}=p_2^{n_2}\cdots p_r^{n_r}$, then clearly, $\mathbb{E}=\mathbb{E}^{\mathcal{G'}}\mathbb{E}^{\mathcal{G}_{p}}$. As $p\nmid [\mathbb{E}:\mathbb{E}^{\mathcal{G'}}]$, in view of (\cite{RS}, Lemma 13.4), the subalgebra $\mathbb{E}*_{\tau_{|_{\mathcal{G}'}}}\mathcal{G'}$ of $\mathbb{E}*_{\tau}\mathcal{G}$ is a split $\mathbb{E}^{\mathcal{G'}}$-algebra, where $\tau_{|_{\mathcal{G'}}}$ is the restriction of $\tau$ on $\mathcal{G'}\times \mathcal{G'}$. As a consequence, in view of (\cite{Her94}, Lemma 4.4.1), the factor set ${\tau_{|_{\mathcal{G'}}}}$ is equivalent to
the trivial factor set. Hence, for each $\sigma_i\in\mathcal{G'}$, replacing $z_{\sigma_i}$ by a suitable $\mathbb{E}$-multiple, we can assume that $\tau_{|_{\mathcal{G'}}}:\mathcal{G'}\times\mathcal{G'}\rightarrow \mathcal{U}(\mathbb{E})$ is the trivial factor set. Therefore, the algebra $\mathbb{E}^{\mathcal{G}_{p}}*_{\tau_{|_{\mathcal{G'}}}}\mathcal{G'}$ make sense, which is clearly split $\mathbb{F}$-algebra and equal to the $\mathbb{F}$-subalgebra of $\operatorname{Cen}_{\mathcal{A}}(\mathcal{B})$ generated by $\mathbb{E}^{\mathcal{G}_{p}}$ and $\mathcal{G'}$. This proves $(i)$.\para\noindent $(ii)$ For the convenience of notation, let us denote $\mathbb{E}^{\mathcal{G}_p}*_{\tau_{|_{\mathcal{G'}}}}\mathcal{G'}$ by $\mathcal{D}$. By (\cite{JdR16}, Theorem 2.1.10), we have that $\operatorname{dim}_{\mathbb{F}}(\operatorname{Cen}_{\operatorname{Cen}_{\mathcal{A}}(\mathcal{B})}(\mathcal{D}))=p^{2n_1}$ and \begin{equation}\label{eqq}
\operatorname{Cen}_{\mathcal{A}}(\mathcal{B})\cong \mathcal{D}\otimes_{\mathbb{F}}\operatorname{Cen}_{\operatorname{Cen}_{\mathcal{A}}(\mathcal{B})}(\mathcal{D}).
\end{equation}Observe that $\mathbb{E}^{\mathcal{G'}}$ commutes with the elements of $\mathbb{E}^{\mathcal{G}_{p}}$ and $\mathcal{G'}$. Therefore, $\mathbb{E}^{\mathcal{G'}}\subseteq \operatorname{Cen}_{\operatorname{Cen}_{\mathcal{A}}(\mathcal{B})}(\mathcal{D})$. As $\operatorname{dim}_{\mathbb{F}}(\operatorname{Cen}_{\operatorname{Cen}_{\mathcal{A}}(\mathcal{B})}(\mathcal{D}))=[\mathbb{E}^{\mathcal{G'}}:\mathbb{F}]^2$, it turns out that $\mathbb{E}^{\mathcal{G'}}$ is a strictly maximal subfield of $\operatorname{Cen}_{\operatorname{Cen}_{\mathcal{A}}(\mathcal{B})}(\mathcal{D})$ and hence $\operatorname{Cen}_{\operatorname{Cen}_{\mathcal{A}}(\mathcal{B})}(\mathcal{D})\cong \mathbb{E}^{\mathcal{G'}}*_{\tau'}\operatorname{Gal}(\mathbb{E}^{\mathcal{G'}}/\mathbb{F})\cong\mathbb{E}^{\mathcal{G'}}*_{\tau'}\mathcal{G}_{p}$ for some factor set $\tau'.$ \para \noindent As the Schur index of $\mathcal{D}$ is 1, we obtain that the Schur index of $\operatorname{Cen}_{\operatorname{Cen}_{\mathcal{A}}(\mathcal{B})}(\mathcal{D})\cong\mathbb{E}^{\mathcal{G'}}*_{\tau'}\mathcal{G}_{p}$ is $p$. This proves $(ii)$.
\para\noindent $(iii)$ This is an immediate consequence of equation (\ref{eqq}). \qed
\para\noindent If $E/F$ is a finite Galois extension, then there exists an element $w\in E$ such that $\{\sigma(w)~|~\sigma \in \G(E/F)\}$ is an $F$-basis of $E$. Such a basis is called {\it a normal basis} and $w$ {\it a normal element} of $E/F$ (see \cite{Lang2012}, Chapter VI, Theorem 13.1).
\begin{theorem}\label{3thm2}
Let $(\mathbb{K}/\mathbb{F},\tau)$ be an abelian crossed product split $\mathbb{F}$-algebra. Let $w$ be a normal element of the extension $\mathbb{K}/\mathbb{F}$. Let $\{ \sigma_{1},\sigma_{2},\ldots,\sigma_{n} \}$ be the Galois group of $\mathbb{K}/\mathbb{F}.$ Then there exist units $\{\mathbf{z}_{1}, \mathbf{z}_{2}, \ldots , \mathbf{z}_{n}\}$ with $\mathbf{z}_{i}\mathbf{z}_{j}=\mathbf{z}_{k}$, which forms a basis of  $(\mathbb{K}/\mathbb{F},\tau)$ as a vector space over $\mathbb{K}$ such that  
	\para \noindent $(i)$ $\{ \mathbf{z}_{i}^{-1} \alpha^{-1} \widehat{E}\, \alpha\, \mathbf{z}_{i}~|~ 1 \leq i \leq n\}$ is a complete set of primitive orthogonal idempotents of $(\mathbb{K}/\mathbb{F},\tau)$. 
	\para \noindent $(ii)$ $\{\mathbf{z}_{i}^{-1} \alpha^{-1} \widehat{E}\, \alpha\, \mathbf{z}_{j}~|~ 1 \leq i,j \leq n\}$ is a complete set of matrix units of $(\mathbb{K}/\mathbb{F},\tau)$.
	\para\noindent Here $\widehat{E}=\frac{1}{n} \sum_{i=1}^{n} \mathbf{z}_{i}$, $\alpha = \sum_{i=1}^{n} \alpha_{i}\mathbf{z}_{i}$ and $\alpha_{1} , \ldots, \alpha_{n} $ in $\mathbb{K}$  satisfy the following relations: 
	\begin{equation*}
	\begin{array}{lll}\noindent
	\sum _{i=1}^{n} \sigma_{i}(w)  \alpha_i&=& \sum_{i=1}^{n}  \sigma_{i}(w)\\
	\sum _{i=1}^{n} ~\sigma_{i}(\sigma_{j}(w) )\alpha_{i}&=& w- \sigma_{j}(w),~ 2 \leq j \leq n.
	\end{array}
	\end{equation*}
\end{theorem}
\noindent \textbf{Proof.} This theorem can be obtained by proceeding similar to the proof of (\cite{BGO}, Lemma 3.2 and Theorem 3.5).~~\qed
\begin{theorem}\label{3thm3}
Let $A=(\mathbb{K}/\mathbb{F},\tau)$ be an abelian crossed product $\mathbb{F}$-algebra such that its dimension (over $\mathbb{F}$) is $p^{2n}$ and the Schur index is $p$. Let $\mathbb{L}$ be a subfield of $\mathbb{K}$ with $[\mathbb{L}:\mathbb{F}]$ least, such that $\mathbb{L}$ splits $A$ and let $[\mathbb{L}:\mathbb{F}]=p^r$ (such an $\mathbb{L}$ exists because $\mathbb{K}$ splits $A$ and $\mathbb{F}$ does not). Then the following statements hold:
\begin{enumerate}[(i)]
\item $(\mathbb{K}/\mathbb{L},\tau_{|_{\operatorname{Gal}(\mathbb{K}/\mathbb{L})}})$ is a split $\mathbb{L}$-subalgebra of $(\mathbb{K}/\mathbb{F}, \tau)$ isomorphic to $M_{p^{n-r}}(\mathbb{L})$, where $\tau_{|_{\operatorname{Gal}(\mathbb{K}/\mathbb{L})}}$ is the restriction of $\tau$ on $\operatorname{Gal}(\mathbb{K}/\mathbb{L})\times\operatorname{Gal}(\mathbb{K}/\mathbb{L})$.
\item If $\{E_{ij}~|~1\leq i,j\leq p^{n-r}\}$ is a complete set of matrix units of $(\mathbb{K}/\mathbb{L},\tau_{|_{\operatorname{Gal}(\mathbb{K}/\mathbb{L})}})$, (which can be computed in view of Theorem \ref{3thm2}) and $B$ is the $\mathbb{F}$-span of $E_{ij}$, $1\leq i,j\leq p^{n-r},$ then $\operatorname{Cen}_{A}(B)$ has the Schur index $p$, its dimension over $\mathbb{F}$ is $p^r$ and $A\cong B\otimes_{\mathbb{F}}\operatorname{Cen}_{A}(B).$
\item $\mathbb{L}$ is a strictly maximal subfield of $\operatorname{Cen}_{A}(B)$ and there exist units $u_{1},u_2,\cdots,u_{p^r}$ with $u_1=1$ such that $\{\sigma_{u_i}~|~1\leq i\leq p^{r}\}$ is the Galois group of $\mathbb{L}/\mathbb{F}$, where $\sigma_{u_i}:\mathbb{L}\rightarrow \mathbb{L}$ is given by $\alpha\mapsto u_{i}^{-1}\alpha u_i$ for all $\alpha\in\mathbb{L}$. Consequently, $\overline{\tau}:\operatorname{Gal}(\mathbb{L}/\mathbb{F})\times \operatorname{Gal}(\mathbb{L}/\mathbb{F})\rightarrow \mathcal{U}(\mathbb{L})$ given by $(u_i,u_j)\mapsto u_k^{-1}u_iu_j,$ where $\sigma_{u_i}\circ\sigma_{u_j}=\sigma_{u_k}$, is a factor set and $\operatorname{Cen}_{A}(B)=\mathbb{L}*_{\bar{\tau}}\operatorname{Gal}(\mathbb{L}/\mathbb{F}).$
\item  Let $\mathbb{L}'$ be a subfield of $\mathbb{L}$ with $[\mathbb{L}:\mathbb{L'}]=p$. Assume that $u_i'$s are ordered so that $\{{\sigma_{u_i}}_{|_{\mathbb{L'}}}~|~1\leq i\leq p^{r-1}\}$ is the Galois group of $\mathbb{L'}/\mathbb{F}$. Then there exist $\beta_i\in \mathbb{L}*_{\bar{\tau}_{|_{\operatorname{Gal}(\mathbb{L}/\mathbb{L}')}}}\operatorname{Gal}(\mathbb{L}/\mathbb{L'}),~1\leq i\leq p^{r-1}$, such that $(\beta_iu_i)(\beta_ju_j)=\beta_ku_k$, if $\sigma_{u_i}\circ \sigma_{u_j}=\sigma_{u_k}.$ Moreover, $\beta_{i}'$s are unique upto multiples of $\mathbb{L'}$. Consequently, $(\mathbb{L'}/\mathbb{F},\tilde{\tau})$ is a split $\mathbb{F}$-subalgebra of $\operatorname{Cen}_{A}(B)$ of dimension $p^{r-1}$ over $\mathbb{F}$, where $\tilde{\tau}:\operatorname{Gal}(\mathbb{L'}/\mathbb{F})\times\operatorname{Gal}(\mathbb{L'}/\mathbb{F})\rightarrow \mathcal{U}(\mathbb{L'})$ is given by $(\beta_iu_i,\beta_ju_j)\mapsto (\beta_ku_k)^{-1}(\beta_iu_i)(\beta_ju_j)$. Here $k$ is such that $\sigma_{u_i}\circ \sigma_{u_j}=\sigma_{u_k}.$ Furthermore, the centralizer of $(\mathbb{L'}/\mathbb{F},\tilde{\tau})$ in $\operatorname{Cen}_{A}(B)$ is a division algebra and hence a complete set of matrix units and primitive orthogonal idempotents of $\operatorname{Cen}_{A}(B)$ is same as those of $(\mathbb{L}'/\mathbb{F},\tilde{\tau})$.
\item Finally, if $\mathcal{E}$ is a complete set of matrix units (respectively primitive orthogonal idempotents) of $(\mathbb{L}'/\mathbb{F},\tilde{\tau})$ obtained using Theorem \ref{3thm2}, then 
\begin{equation}
\{eE_{ij}~|~e\in\mathcal{E},1\leq i,j\leq p^{n-r}\} ~~~ (\mbox{respectively}~\{eE_{ii}~|~e\in\mathcal{E},1\leq i\leq p^{n-r}\})
\end{equation}
is a complete set of matrix units (respectively primitive orthogonal idempotents) of $(\mathbb{K}/\mathbb{F}, \tau)$.
\end{enumerate}
\end{theorem}
\noindent \textbf{Proof.} $(i)$ Observe that \begin{equation}\label{eq1}
(\mathbb{K}/\mathbb{L},\tau_{|_{\operatorname{Gal}(\mathbb{K}/\mathbb{L})}})\subseteq\operatorname{Cen}_{A}(\mathbb{L}).
\end{equation}  By using (\cite{JdR16}, Proposition 2.2.4 (1)), we have $\operatorname{dim}_{\mathbb{L}}(\operatorname{Cen}_{A}(\mathbb{L}))=\frac{\operatorname{dim}_{\mathbb{F}}(A)}{[\mathbb{L}:\mathbb{F}]^2}=\frac{p^{2n}}{p^{2r}}=p^{2(n-r)}.$ Also, $\operatorname{dim}_{\mathbb{L}}(\mathbb{K}/\mathbb{L},\tau_{|_{\operatorname{Gal}(\mathbb{K}/\mathbb{L})}})=|\operatorname{Gal}(\mathbb{K}/\mathbb{L})|^2=p^{2(n-r)}.$ Therefore, the equality holds in equation (\ref{eq1}), i.e., \begin{equation}\label{eqqq}
(\mathbb{K}/\mathbb{L},\tau_{|_{\operatorname{Gal}(\mathbb{K}/\mathbb{L})}})=\operatorname{Cen}_{A}(\mathbb{L}).
\end{equation} By (\cite{JdR16}, Proposition 2.2.4 (3)), $\operatorname{Cen}_{A}(\mathbb{L})$ is Brauer equivalent to $A\otimes_{\mathbb{F}}\mathbb{L},$ which splits, as $\mathbb{L}$ is a splitting field of $A$. Consequently, $(\mathbb{K}/\mathbb{L},\tau_{|_{\operatorname{Gal}(\mathbb{K}/\mathbb{L})}})$ is a split $\mathbb{L}$-subalgebra of $(\mathbb{K}/\mathbb{F},\tau)$. As the $\operatorname{dim}_{\mathbb{L}}(\mathbb{K}/\mathbb{L},\tau_{|_{\operatorname{Gal}(\mathbb{K}/\mathbb{L})}})=p^{2(n-r)}$, we have that $(\mathbb{K}/\mathbb{L},\tau_{|_{\operatorname{Gal}(\mathbb{K}/\mathbb{L})}})$ is isomorphic to $M_{p^{n-r}}(\mathbb{L})$. This proves $(i).$
\para \noindent $(ii)$ By using (\cite{JdR16}, Lemma 2.6.2),
$$A\cong (\mathbb{K}/\mathbb{L},\tau_{|_{\operatorname{Gal}(\mathbb{K}/\mathbb{L})}})*\operatorname{Gal}(\mathbb{K}/\mathbb{F})/\operatorname{Gal}(\mathbb{K}/\mathbb{L}).$$\para \noindent If $\{E_{ij}~|~1\leq i,j\leq p^{n-r}\}$ is a complete set of matrix units of $(\mathbb{K}/\mathbb{L},\tau_{|_{\operatorname{Gal}(\mathbb{K}/\mathbb{L})}})$ and $B$ is the $\mathbb{F}$-span of $E_{ij},~1\leq i,j\leq p^{n-r}$, then $B$ is a $\mathbb{F}$-subalgebra of $(\mathbb{K}/\mathbb{L},\tau_{|_{\operatorname{Gal}(\mathbb{K}/\mathbb{L})}})$ isomorphic to $M_{p^{n-r}}(\mathbb{F})$. By using (\cite{JdR16}, Theorem 2.1.10), we have that \begin{equation}\label{eq2}
A\cong B\otimes_{\mathbb{F}}\operatorname{Cen}_{A}(B)
\end{equation} and $\operatorname{dim}_{\mathbb{F}}(\operatorname{Cen}_{A}(B))=\frac{\operatorname{dim}_{\mathbb{F}}(A)}{\operatorname{dim}_{\mathbb{F}}(B)}=\frac{p^{2n}}{p^{2(n-r)}}=p^{2r}.$ Since the Schur index of $B$ is $1$, it turns out from equation (\ref{eq2}) that the Schur index of $\operatorname{Cen}_{A}(B)$ is $p$. This proves $(ii)$. \para \noindent $(iii)$ Observe that $\mathbb{L}\subseteq\operatorname{Cen}_{A}(B),$ as the scalar matrices with all entries in $\mathbb{L}$ commute with the elements of $B.$ Since $[\mathbb{L}:\mathbb{F}]$ is $p^r$ and $\operatorname{dim}_{\mathbb{F}}(\operatorname{Cen}_{A}(B))$ is $p^{2r}$, it turns out that $\mathbb{L}$ is a strictly maximal subfield of $\operatorname{Cen}_{A}(B).$ Suppose $\operatorname{Gal}(\mathbb{K}/\mathbb{F})=\{\sigma_{1},\cdots,\sigma_{p^n}\}$ and $\sigma_i$'s are ordered so that $\{{\sigma_{i}}_{|_{\mathbb{L}}}$, $1\leq i\leq p^{r}\}$ is the Galois group of $\mathbb{L}/\mathbb{F}$. By (\cite{JdR16}, Theorem 2.1.9), for each $\sigma_{i}{_{|_{\mathbb{L}}}} \in \operatorname{Gal}(\mathbb{L}/\mathbb{F}),~1\leq i\leq p^{r},$ there is a unit $u_i$ in $\operatorname{Cen}_{A}(B)$ such that ${\sigma_{i}}_{|_{\mathbb{L}}}=\sigma_{u_i}$, i.e., the conjugation automorphism on $\mathbb{L}$ by $u_{i}$. As $\operatorname{dim}_{\mathbb{L}}(\C_{A}(B))= |\G(\mathbb{L}/\mathbb{F})|=p^r,$ it follows that $\{u_{i} ~|~ \sigma_i \in \G(\mathbb{L}/\mathbb{F})\}$ is a basis for the $\mathbb{L}$-vector space $\C_{A}(B)$. If $\bar{\tau}$ is as given in the statement, then the associativity of multiplication in $\operatorname{Cen}_{A}(B)$ implies that $\bar{\tau}$ is a factor set and  $$\operatorname{Cen}_{A}(B)\cong \mathbb{L}*_{\overline{\tau}}\operatorname{Gal}(\mathbb{L}/\mathbb{F}).$$ \noindent $(iv)$ As the Schur index of $\operatorname{Cen}_{A}(B)$ is $p$, we have that $\operatorname{Cen}_{A}(B)\cong M_{p^{r-1}}(D)$ for some $\mathbb{F}$-division algebra $D$ with $\operatorname{dim}_{\mathbb{F}}(D)=p^2$. This gives $\operatorname{Cen}_{A}(B)\otimes_{\mathbb{F}}\mathbb{L}'\cong M_{p^{r-1}}(D\otimes_{\mathbb{F}}\mathbb{L'}).$ By (\cite{JdR16}, Proposition 2.2.4 (3)), $\operatorname{Cen}_{A}(B)\otimes_{\mathbb{F}}\mathbb{L'}$ is Brauer equivalent to $\operatorname{Cen}_{\operatorname{Cen}_{A}(B)}(\mathbb{L'})$, which can be seen to be equal to $(\mathbb{L}/\mathbb{L'},\bar{\tau}_{|_{\operatorname{Gal}(\mathbb{L}/\mathbb{L'})}})$ using the argument similar to those given in equation (\ref{eqqq}). Now if $\mathbb{L'}$ splits $\operatorname{Cen}_{A}(B)$, then in view of equation (\ref{eq2}), $\mathbb{L'}$ also splits $A$, which is not so. Therefore, $\mathbb{L'}$ does not split $\operatorname{Cen}_{A}(B).$ Consequently, $\operatorname{Cen}_{A}(B)\otimes_{\mathbb{F}}\mathbb{L'}$ and hence $(\mathbb{L}/\mathbb{L'},\bar{\tau}_{|_{\operatorname{Gal}(\mathbb{L}/\mathbb{L'})}})$ is not a split $\mathbb{L'}$-algebra. But $\operatorname{dim}_{\mathbb{L'}}(\mathbb{L}/\mathbb{L'},\bar{\tau}_{|_{\operatorname{Gal}(\mathbb{L}/\mathbb{L'})}})=p^2.$ This gives $(\mathbb{L}/\mathbb{L'},\bar{\tau}_{|_{\operatorname{Gal}(\mathbb{L}/\mathbb{L'})}})$ is a division $\mathbb{L'}$- algebra. For notational convenience, let us denote $(\mathbb{L}/\mathbb{L'},\bar{\tau}_{|_{\operatorname{Gal}(\mathbb{L}/\mathbb{L'})}})$ by $D'$. Since $D'$ and $D\otimes_{\mathbb{F}}\mathbb{L'}$ are Brauer equivalent and $\operatorname{dim}_{\mathbb{L'}}(D')=\operatorname{dim}_{\mathbb{L'}}(D\otimes_{\mathbb{F}}\mathbb{L'})=p^2,$ in view of (\cite{RS}, Proposition 12.5b (i)), $D'$ and $D\otimes_{\mathbb{F}}\mathbb{L'}$ are isomorphic as $\mathbb{L'}$-algebras. Let $\phi:D\otimes_{\mathbb{F}}\mathbb{L'}\rightarrow D'$ be an $\mathbb{L'}$-algebra isomorphism. Now, as $D$ is a division algebra with center $\mathbb{F}$, $\phi(D\otimes_{\mathbb{F}}\mathbb{F})$ is also a division algebra with center $\mathbb{F}$ and is inside $D'$. Let us denote $\phi(D\otimes_{\mathbb{F}}\mathbb{F})$ by $\mathbb{D}.$ As $D\cong \mathbb{D}$ as $\mathbb{F}$-algebras, we have that
\begin{equation}\label{eq11}
\operatorname{Cen}_{A}(B)\cong M_{p^{r-1}}(D)\cong M_{p^{r-1}}(\mathbb{D})~\mbox{as}~\mathbb{F}\mbox{-algebras}.
\end{equation} By (\cite{JdR16}, Theorem 2.1.10), we have that \begin{equation}\label{eq12}
\operatorname{Cen}_{A}(B)\cong \mathbb{D}\otimes_{\mathbb{F}}\operatorname{Cen}_{\operatorname{Cen}_{A}(B)}(\mathbb{D})
\end{equation} and a complete set of matrix units (respectively primitive orthogonal idempotents) of $\operatorname{Cen}_{A}(B)$ is precisely a complete set of matrix units (respectively primitive orthogonal idempotents) of $\operatorname{Cen}_{\operatorname{Cen}_{A}(B)}(\mathbb{D})$. \para \noindent We now assert that $\operatorname{Cen}_{\operatorname{Cen}_{A}(B)}(\mathbb{D})$ is a split $\mathbb{F}$-algebra isomorphic to $(\mathbb{L'}/\mathbb{F}, \tilde{\tau})$, where $\tilde{\tau}$ is as given in statement. Firstly, we see that $\operatorname{Cen}_{\operatorname{Cen}_{A}(B)}(\mathbb{D})\cong M_{p^{r-1}}(\mathbb{F})$. From equations (\ref{eq11}) and (\ref{eq12}), we have that $\mathbb{D}\otimes_{\mathbb{F}} \operatorname{Cen}_{\operatorname{Cen}_{A}(B)}(\mathbb{D}) \cong M_{p^{r-1}}(\mathbb{D})$. This gives $\mathbb{D^{\circ}}\otimes_{\mathbb{F}}\mathbb{D}\otimes_{\mathbb{F}}\operatorname{Cen}_{\operatorname{Cen}_{A}(B)}(\mathbb{D})\cong \mathbb{D^{\circ}}\otimes_{\mathbb{F}} M_{p^{r-1}}(\mathbb{D})$, where $\mathbb{D}^{\circ}$ is the opposite algebra of $\mathbb{D}$. Hence, in view of (\cite{JdR16}, Problem 2.4.1), we obtain that the Schur index of $\operatorname{Cen}_{\operatorname{Cen}_{A}(B)}(\mathbb{D})$ is $1$ and hence $\operatorname{Cen}_{\operatorname{Cen}_{A}(B)}(\mathbb{D})\cong M_{p^{r-1}}(\mathbb{F})$.\para \noindent Finally, we will show that $\operatorname{Cen}_{\operatorname{Cen}_{A}(B)}(\mathbb{D})\cong (\mathbb{L'}/\mathbb{F}, \tilde{\tau})$, where $\tilde{\tau}$ is as given in the statement. As $\mathbb{D}\subset D'$ and the center of $D'$ is $\mathbb{L'}$, it follows that $\mathbb{L'}\subset \operatorname{Cen}_{\operatorname{Cen}_{A}(B)}(\mathbb{D}).$ Now $\operatorname{dim}_{\mathbb{F}}(\operatorname{Cen}_{\operatorname{Cen}_{A}(B)}(\mathbb{D}))=p^{2(r-1)}=[\mathbb{L'}:\mathbb{F}]^2,$ it turns out that $\mathbb{L'}$ is a strictly maximal subfield of $\operatorname{Cen}_{\operatorname{Cen}_{A}(B)}(\mathbb{D}).$ Let us reorder $u_i$$'$s, $1\leq i\leq p^{r},$ so that ${\sigma_{u_i}}_{|_\mathbb{L'}},~1\leq i\leq p^{r-1}$, is the Galois group of $\mathbb{L'}/\mathbb{F}.$ For each ${\sigma_{u_i}}_{|_\mathbb{L'}}\in\operatorname{Gal}(\mathbb{L'}/\mathbb{F}),$ by (\cite{JdR16}, Theorem 2.1.9), there exists $w_i\in \operatorname{Cen}_{\operatorname{Cen}_{A}(B)}(\mathbb{D})$ such that ${\sigma_{u_i}}_{|_{\mathbb{L'}}}(\alpha)=w_{i}^{-1}\alpha w_{i},$ for all $\alpha\in\mathbb{L'}$. Hence, $\operatorname{Cen}_{\operatorname{Cen}_{A}(B)}(\mathbb{D})\cong(\mathbb{L'}/\mathbb{F},\tilde{\tau})$, where $\tilde{\tau}:\operatorname{Gal}(\mathbb{L'}/\mathbb{F})\times \operatorname{Gal}(\mathbb{L'}/\mathbb{F})\rightarrow \mathcal{U}(\mathbb{L'})$ is given by $(w_i,w_j) \mapsto w_{k}^{-1}w_iw_j$. Since $\operatorname{Cen}_{\operatorname{Cen}_{A}(B)}(\mathbb{D})$ is a split $\mathbb{F}$-algebra, the factor set $\tilde{\tau}$ is equivalent to the trivial factor set. Hence, by replacing $w_i$ with a suitable multiple of $\mathbb{L'}$, we may assume that $w_iw_j=w_k$, whenever ${\sigma_{u_i}}_{|_{\mathbb{L'}}}\circ {\sigma_{u_j}}_{|_{\mathbb{L'}}}={\sigma_{u_k}}_{|_{\mathbb{L'}}}$, which is equivalent to saying that $\sigma_{u_i}\circ\sigma_{u_j}=\sigma_{u_k}.$ Further, as ${\sigma_{u_i}}_{|_{\mathbb{L'}}}(\alpha)=w_i^{-1}\alpha w_i$, we obtain that for all $1\leq i\leq p^{r-1}$, $w_iu_i^{-1}$ belongs to $\operatorname{Cen}_{\operatorname{Cen}_{A}(B)}(\mathbb{L'})$, which is equal to $(\mathbb{L}/\mathbb{L'},\tilde{\tau}_{|_{\operatorname{Gal}(\mathbb{L}/\mathbb{L'})}})$. Consequently, $w_i=\beta_iu_i$ for some $\beta_i\in(\mathbb{L}/\mathbb{L'},\tilde{\tau}_{|_{\operatorname{Gal}(\mathbb{L}/\mathbb{L'})}}).$ For uniqueness of $\beta_i$$'$s up to a multiple of $\mathbb{L'}$, let us assume that $\beta_i$ and $\beta_i'$ belonging to $(\mathbb{L}/\mathbb{L'},\tilde{\tau}_{|_{\operatorname{Gal}(\mathbb{L}/\mathbb{L'})}})$ are such that $(\beta_iu_i)^{-1}\alpha \beta_iu_i=(\beta_i'u_i)^{-1}\alpha \beta_i'u_i$ for all $\alpha \in \mathbb{L'}$. This gives $\beta_i'u_i(\beta_iu_i)^{-1}\alpha \beta_iu_i (\beta_i'u_i)^{-1}=\alpha$ for all $\alpha \in \mathbb{L'}$ and hence $\beta_i'\beta_i^{-1}$ belongs to the centralizer of $\mathbb{L'}$ in $\operatorname{Cen}_{\operatorname{Cen}_{A}(B)}(\mathbb{D})$, which is $\mathbb{L'}$, as it is a strictly maximal subfield of $\operatorname{Cen}_{\operatorname{Cen}_{A}(B)}(\mathbb{D})$. This proves $(iv)$.\para\noindent $(v)$ In $(ii),$ we have shown that $A\cong B\otimes_{\mathbb{F}} \operatorname{Cen}_{A}(B)$ and also $\{E_{ij}~|~1\leq i,j\leq p^{n-r}\}$ is a complete set of matrix units of $B$. In $(iv)$, we have shown that a complete set of matrix units and primitive orthogonal idempotents of $\operatorname{Cen}_{A}(B)$ is precisely same as those of $(\mathbb{L'}/\mathbb{F},\tilde{\tau}).$ Therefore, $(v)$ follows.~~\qed
\section{Some Examples}
In this section, we will provide some examples of groups $G$ for which a complete set of matrix units and, in particular, primitive orthogonal idempotents of $\Q G$ can be computed using the method developed in section 3.
\subsection{A class of Frobenius groups}
In (\cite{BGO}, Theorem 6.1), we proved that a Frobenius group $G$ of odd order with cyclic complement is a generalized strongly monomial group and the Schur indices of all the simple components of $\Q G$ is 1. Using essentially the same idea with obvious modifications, we prove the following:
\begin{theorem}
If $G$ is a monomial Frobenius group with complement $K$ such that the Schur indices of all the simple components of $\Q K$ is either $1$ or a prime, then $G$ is a generalized strongly monomial group and the Schur indices of all the simple components of $\Q G$ is either $1$ or a prime.
\end{theorem} 
\noindent \textbf{Proof.} Let $G$ be a Frobenius group with kernel $N$ and complement $K$. From (\cite{Hupp98}, Theorems 16.7 a) and 18.7), it is known that $N$ is nilpotent. If $\chi$ is an irreducible character of $G$, then $\chi$ is either lifted from an irreducible character, say $\psi$, of $K\cong G/N$ or it is induced from an irreducible character, say $\phi$, of $N$. If $\chi$ is lifted from $\psi$, then $\Q Ge_{\Q}(\chi)\cong \Q Ke_{\Q}(\psi)$. As the Schur index of $\Q Ke_{\Q}(\psi)$ is either $1$ or a prime, the same holds for $\Q Ge_{\Q}(\chi)$. If $\chi$ is induced from $\phi$, then arguing as in (\cite{BGO}, Theorem 6.1), we have that the Schur index of $\Q Ge_{\Q}(\chi)$ is same as that of $\Q Ne_{\Q}(\phi)$, which is either $1$ or $2$ as $N$ is nilpotent. Hence, the result is proved.~~\qed
\begin{cor}
A complete set of matrix units and, in particular, primitive orthogonal idempotents of $\Q G$ can be computed using the method developed in section 3, if $G$ is one of the following kind:
\begin{enumerate}[(i)]
\item Nilpotent-by-Nilpotent group with fixed-point-free action.
\item Camina groups.
\end{enumerate}
\end{cor}
\noindent \textbf{Proof.} $(i)$ Suppose $G=N\rtimes K$, where $N$ and $K$ are nilpotent groups and the action of $K$ on $N$ is fixed-point-free. Such a $G$ is a Frobenius group (see \cite{Hupp98}, Proposition 16.5). Also, it is monomial because a Frobenius group is monomial if and only if its complement is monomial. Further, the Schur indices of all the simple components of $\Q K$ is either $1$ or $2$ as $K$ is nilpotent. Therefore, the result follows from the above theorem.\para\noindent $(ii)$ If $G$ is a Camina group, then either of the following holds (see \cite{lewis2014}): 
\begin{enumerate}[(a)]
\item $G$ is a $p$-group;
\item $G$ is a Frobenius group whose complement is either cyclic or is isomorphic to the group of quaternions.
\end{enumerate}  If $G$ is of type (a), then all the simple components of $\Q G$ have the Schur indices either $1$ or $2$. If $G$ is of type (b), then the complement $K$ of $G$ being cyclic or quaternion, the Schur indices of all the simple components of $\Q K$ is either $1$ or $2$. Therefore, the result follows from the above theorem.~~\qed
\subsection{Groups of order $p_1^2p_2^2$ and $p_1^{\alpha}p_2\cdots p_k,$ where $\alpha=1,2$}
\begin{theorem}
A complete set of matrix units and, in particular, primitive orthogonal idempotents of $\Q G$ can be computed using the method developed in Section 3, if $G$ is one of the following kind:
\begin{enumerate}[(i)]
\item Groups of order $p_1p_2\cdots p_k$ and $p_1^2p_2\cdots p_k$, where $p_{i}$$'$s are distinct primes.
\item Groups of order $p_{1}^{2}p_{2}^2$.
\end{enumerate}
\end{theorem}
\noindent \textbf{Proof.} $(i)$ Let $G$ be a group of order $p_1^{\alpha}p_2\cdots p_k,$ where $\alpha$ is either $1$ or $2$. Let $\chi\in \operatorname{Irr}(G)$ and let $m_{\Q}(\chi)$ be the Schur index of the simple component $\Q Ge_{\Q}(\chi)$ of $\Q G$. For $i\geq 2,$ in view of (\cite{Isa76}, Theorem 10.9), $p_i$ cannot divide $m_{\Q}(\chi)$ because a Sylow $p_i$-subgroup of $G$ is cyclic group of order $p_i$ and hence an elementary abelian group. Also, a Sylow $p_1$-subgroup of $G$ is one of the following: $C_{p_1}$, cyclic of order $p_1$ (which is the case when $\alpha=1$) or $C_{p_1}\times C_{p_1}$ or $C_{p_1^2}$. In the first two cases, again by (\cite{Isa76}, Theorem 10.9), $p_1$ cannot divide $m_{\Q}(\chi)$. In the last case, when a Sylow $p_1$-subgroup of $G$ is $C_{p_1^2}$, we have that if $p_1$ divides $m_{\Q}(\chi)$, then $p_1m_{\Q}(\chi)$ must divide $|G|$. Consequently, $m_{\Q}(\chi)=1$ or $p_1$.
\para\noindent $(ii)$ Let $G$ be a group of order $p_{1}^2p_{2}^2.$ If $p_1=p_2$, then $G$ being a $p_1$-group, the Schur indices of all the simple components of $\Q G$ is either $1$ or $2$. If $p_1\neq p_2$, then without loss of generality, we may assume that $p_{1}>p_{2}.$ Let $G_{p_{i}}$ be a Sylow $p_{i}$-subgroup of $G$. \para\noindent For $i=1,2$, if $G_{p_i}=C_{p_i}\times C_{p_i},$ then $p_i$ cannot divide $m_{\Q}(\chi)$ for any $\chi\in \operatorname{Irr}(G)$. Hence, $m_{\Q}(\chi)=1.$ \para \noindent If $G_{p_1}=C_{p_1^2}$ and $G_{p_2}=C_{p_2}\times C_{p_2},$ then $p_2$ cannot divide $m_{\Q}(\chi)$ and if $p_1$ divide $m_{\Q}(\chi)$, then we have that $p_1m_{\Q}(\chi)$ must divide $|G|.$ This gives $m_{\Q}(\chi)=1$ or $p_1.$ Similarly, if $G_{p_1}=C_{p_1}\times C_{p_1}$ and $G_{p_2}=C_{p_2^2},$ then $m_{\Q}(\chi)=1$ or $p_2.$\para\noindent Finally, suppose $G_{p_i}=C_{p_i^2}$ for $i=1,2$. If $(p_1,p_2)= (3,2)$, by using GAP, one can see that the Schur index of any simple component of $\Q G$ is either $1$ or $2$. Let us suppose $(p_1,p_2)\neq (3,2)$. In this case, one can show that $G_{p_1}$ is the normal subgroup of $G$ and hence $G$ has the following presentation: $$G=\langle a,b~|~ a^{p_1^2}=b^{p_2^2}=1,b^{-1}ab=a^r\rangle, $$ where $\operatorname{gcd}(r,p_1)=1$ and $p_1^2|r^{p_2^2}-1$. From (\cite{JdR16}, Theorem 3.5.12), we know that any $\chi\in \operatorname{Irr}(G)$ is arising from a Shoda pair of the type $(H,K)$, where $H=\langle a,b^d\rangle $, $d$ divides the order of $r$ modulo $p_1^2$. Also, from (\cite{JdR16}, Theorem 3.5.5), $\Q Ge_{\Q}(\chi)$ is isomorphic to $\Q He_{\Q}(\lambda)*N_{G}(K)/H$ and hence $m_{\Q}(\chi)$ divides $|N_{G}(K)/H|$. As $\langle a\rangle \subseteq H,$ we have that $|N_{G}(K)/H|$ is a divisor of $p_2^2$. Thus $p_1$ cannot divide $m_{\Q}(\chi)$ and if $p_2$ divide $m_{\Q}(\chi)$, then we have that $p_2m_{\Q}(\chi)$ must divide $|G|.$ This implies $m_{\Q}(\chi)=1$ or $p_2$. ~~\qed
\subsection{$G=((C_p\times C_p)\rtimes C_p)\rtimes C_{2^n}$}
In (\cite{BG}, Section 4), we considered a group $G$ of the type $((C_p\times C_p)\rtimes C_p)\rtimes C_{2^n}$, where $n\geq 2$, $p$ is an odd prime and $2^{n-1}$ is the highest power of 2 dividing $p-1$. For the precise presentation of $G$, we refer to section 4 of \cite{BG}. In that work, it is proved that $G$ is a generalized strongly monomial group, all the simple components of $\Q G$ have the Schur indices 1 except one component, where the Schur index is 2. Thus, a complete set of matrix units and, in particular, primitive orthogonal idempotents of $\Q G$ can be computed using the method developed in section 3. \para\noindent Let us now illustrate the application of the method developed in Section 3 when $p=5$ (in which case $G$ is SmallGroup(1000,86) in GAP library).\para\noindent When $p=5$, $G$ is generated by $x,y,z,w$ with the following defining relations:
\begin{equation*}
\begin{array}{lll}
x^5=y^5=z^5=w^8=1,\\
xy=yx,xz=zx,yz=zyx,\\
w^{-1}xw=x^2, w^{-1}yw=z^3,w^{-1}zw=y.
\end{array}
\end{equation*}
A complete set of matrix units and primitive orthogonal idempotents of all the simple components of $\Q G$ having the Schur indices $1$ can be computed using earlier work in \cite{BGO}. From (\cite{BG}, Section 4, Lemma 4), the simple component of $\Q G$ having the Schur index 2 is arising from the generalized strong Shoda pair $(H,K)$, where  $H=\langle x,y,w^{4}\rangle$, $K=\langle y\rangle$ and this simple component is isomorphic to $M_{5}(\Q(\zeta_5)/\Q,\sigma,-1) (\cong M_{10}(H(\Q)))$, which is the algebra  $\mathcal{A}$ under consideration. Here $\zeta_5$ is a $5$-th root of unity. From the definition of $\mathcal{B}$, we have that $\mathcal{B}\cong M_5(\Q)$ and $\operatorname{Cen}_{\mathcal{A}}(\mathcal{B})\cong (\Q(\zeta_5)/\Q,\sigma,-1)$.\para\noindent
\textbf{Step I:} A complete set of matrix units and primitive orthogonal idempotents of $\mathcal{B}$.\\
For this part, we will apply Theorem \ref{l1}. The first requirement is that of a strong inductive chain from $H$ to $G$. In (\cite{BG}, Lemma 4), it is proved that $H_0=H\leq H_1=\langle x,y,z,w^4\rangle \leq G$ is a strong inductive chain from $H$ to $G$ with $C_0=H_0$ and $C_1=G$. The next requirement is that of a left transversal set $T_0$ of $C_0$ in $H_1$ and $T_1$ of $C_1$ in $G$. Clearly, we can take $T_0=\{z^i~|~1\leq i\leq 5\}$ and $T_1=\{1\}$. Consequently, the following is obtained:\para\noindent Matrix units of $\mathcal{B}$: $\{z^{-i}\varepsilon(H,K)z^j~|~1\leq i,j\leq 5\}.$\para\noindent Primitive orthogonal idempotents of $\mathcal{B}$:  $\{z^{-i}\varepsilon(H,K)z^i~|~1\leq i\leq 5\}.$ \para \noindent
\textbf{Step II:} A complete set of matrix units and primitive orthogonal idempotents of $\operatorname{Cen}_{\mathcal{A}}(\mathcal{B})$.\\ The dimension of $\operatorname{Cen}_{\mathcal{A}}(\mathcal{B})$ over $\mathbb{F}$ is $16$ and the Schur index is $2$. Therefore, we can directly apply Theorem \ref{3thm3}, when $A=(\Q(\zeta_5)/\Q,\sigma,-1)$ and $p=2$. In this case, $\mathbb{K}=\Q(\zeta_5)$ and $\mathbb{F}=\Q.$ The only non-trivial subfield of $\mathbb{K}$ containing $\mathbb{F}$ is $\Q(\zeta_5+\zeta_5^{-1}),$ which does not split $\operatorname{Cen}_{\mathcal{A}}(\mathcal{B})$ because $-1\notin N_{\Q(\zeta_5)/\Q(\zeta_5+\zeta_5^{-1})}(\Q(\zeta_5) ^{\times})$ (see \cite{BG}, Step 6). Hence, the only choice of $\mathbb{L}$ is $\mathbb{K}$, $[\mathbb{L}:\mathbb{F}]=4$, and $r=2.$ The computations done for different parts of Theorem \ref{3thm3} are as follows:
\para \noindent $(i)$  The subalgebra $(\mathbb{K}/\mathbb{L},\tau_{|_{\operatorname{Gal}(\mathbb{K}/\mathbb{L})}})$ is $\mathbb{L},$ as $\mathbb{L}=\mathbb{K}.$
\para \noindent $(ii)$ From $(i)$, $\{\varepsilon(H,K)\}$ is a complete set of matrix units of $(\mathbb{K}/\mathbb{L},\tau_{|_{\operatorname{Gal}(\mathbb{K}/\mathbb{L})}})$. Therefore, $B=\mathbb{F}$ and hence $\operatorname{Cen}_{A}(B)=A=(\Q(\zeta_5)/\Q,\sigma,-1)$. The Schur index of $\operatorname{Cen}_{A}(B)$ is 2, its dimension over $\mathbb{F}$ is $16$, and $A=\mathbb{F}\otimes_{\mathbb{F}}\operatorname{Cen}_{A}(B).$\para\noindent $(iii)$ 
In order to understand the action and the factor set of the cyclotomic algebra $\operatorname{Cen}_{A}(B)$, the aim of this part is to find a set of units $\{u_1,u_2,u_3,u_4\}$ of $\operatorname{Cen}_{A}(B)$, which forms a basis of $\operatorname{Cen}_{A}(B)$ as a vector space over $\Q (\zeta_5)$. Such a set of $u_i$$'$s has already been computed in (\cite{BG}, Lemma 4) and we recall that the following set of $u_i$$'$s work:$$u_i=\left(\frac{z_w}{\sqrt{5}}\right)^{i},~0\leq i\leq 3,$$ where $z_w=wA_w$ and  $$A_w=\begin{pmatrix}
\varepsilon & \varepsilon & \varepsilon & \varepsilon & \varepsilon\\
\varepsilon & x^4\varepsilon & x^3\varepsilon & x^2\varepsilon & x\varepsilon\\
\varepsilon & x^3\varepsilon & x\varepsilon & x^4\varepsilon & x^2\varepsilon\\
\varepsilon & x^2\varepsilon & x^4\varepsilon & x\varepsilon & x^3\varepsilon\\
\varepsilon & x\varepsilon & x^2\varepsilon & x^3\varepsilon & x^4\varepsilon
\end{pmatrix}$$ with $\varepsilon=\frac{\varepsilon(H,K)}{5}$.  
\para\noindent $(iv)$ Clearly $\mathbb{L'}=\Q(\zeta_5+\zeta_5^{-1})$. Also, the conjugation by $u_0$ and $u_1$ on $\mathbb{L'}$ gives the Galois group of $\mathbb{L'}/\mathbb{F}.$ The next hurdle is to find $\beta_0,\beta_1$ in $(\mathbb{L}/\mathbb{L'},\tau_{|_{\operatorname{Gal}(\mathbb{L}/\mathbb{L'})}})$ such that $(\beta_1u_1)^2=\beta_0u_0=1.$ Using some computations, we see that $\beta_0=1$ and $\beta_1=\frac{(\zeta_5-\zeta_5^{-1})(u_2+(\zeta_5+\zeta_5^{-1}))}{\sqrt{5}}$ works. Hence, $(\mathbb{L'}/\mathbb{F},\tilde{\tau})=(\mathbb{L'}/\mathbb{F},\beta_1u_1,1)$ is a split $\mathbb{F}$-subalgebra of $\operatorname{Cen}_{A}(B)$ and it is of dimension $4$ over $\mathbb{F}$.
\para \noindent $(v)$ In view of $(i)$, a complete set of matrix units and primitive orthogonal idempotents of $\operatorname{Cen}_{A}(B)$ is same as those of $(\mathbb{L}/\mathbb{F},\tilde{\tau}),$ which can be computed using Theorem \ref{3thm2}. Let us illustrate the same here:\para\noindent As $\{\beta_0u_0,\beta_1u_1\}$ is a set of units of $\operatorname{Cen}_{A}(B)$ (as a vector space over $\Q(\zeta_5)$), we can take $\mathbf{z}_1=\beta_0u_0=1$, $\mathbf{z}_2=\beta_1u_1$ and $\widehat{E}=\frac{1+\beta_1u_1}{2}$.\para\noindent Note that $[\mathbb{L'}:\mathbb{F}]=2$ and $\{(1+\sqrt{5}),(\beta_1u_1)^{-1}(1+\sqrt{5})(\beta_1u_1)\}$ is a linearly independent set over $\mathbb{F}$ because $(\beta_1u_1)^{-1}(1+\sqrt{5})(\beta_1u_1)=1-\sqrt{5}$. Hence, $w=1+\sqrt{5}$ is a normal element of the extension $\mathbb{L'}/\mathbb{F}.$\para\noindent Another ingredient required to apply Theorem \ref{3thm2} is $\{\alpha_1,\alpha_2 \}$ in $\mathbb{L'}$ such that the following system of equations hold: $$\begin{pmatrix}
1+\sqrt{5} & 1-\sqrt{5}\\
1-\sqrt{5} & 1+\sqrt{5}
\end{pmatrix}\begin{pmatrix}
\alpha_1\\
\alpha_2
\end{pmatrix}=\begin{pmatrix}
2\\
2\sqrt{5}
\end{pmatrix}.$$ On solving, we get $\alpha_1=\frac{3}{\sqrt{5}}$ and $\alpha_2=\frac{2+\sqrt{5}}{\sqrt{5}}.$ Consequently, if $\alpha=\alpha_1\mathbf{z}_1+\alpha_2\mathbf{z}_2$, then the following holds:\para\noindent Matrix units of $(\mathbb{L'}/\mathbb{F},\tilde{\tau})$: $\{\alpha^{-1}\widehat{E}\alpha,(\beta_1u_1)^{-1}\alpha^{-1}\widehat{E}\alpha,\alpha^{-1}\widehat{E}\alpha \beta_1u_1,(\beta_1u_1)^{-1}\alpha^{-1}\widehat{E}\alpha\beta_1u_1\}.$ \para\noindent Primitive orthogonal idempotents of  $(\mathbb{L'}/\mathbb{F},\tilde{\tau})$: $\{\alpha^{-1}\widehat{E}\alpha,(\beta_1u_1)^{-1}\alpha^{-1}\widehat{E}\alpha\beta_1u_1\}.$ 

\para \noindent \textbf{Step III:} Finally, a complete set of matrix units (respectively primitive orthogonal idempotents) of $\mathcal{A}$.\\
Using steps I and II$(v)$ and using the fact that $\mathcal{A}\cong \mathcal{B}\otimes_{\mathbb{F}}\operatorname{Cen}_{\mathcal{A}}(\mathcal{B})$, we obtain the following:\para \noindent
A complete set of matrix units of $\mathcal{A}$ is $$\{z^{-i}\varepsilon(H,K)z^{j}\varepsilon(H,K)(\beta_1u_1)^{-i'}\alpha^{-1}\widehat{E}\alpha(\beta_1u_1)^{j'}~|~1\leq i,j\leq 5,1\leq i',j'\leq 2\}.$$
By using (\cite{BK19}, Lemma 4), we have that $e_{\Q}(\lambda^G)\varepsilon(H,K)=\varepsilon(H,K)=\varepsilon(H,K)e_{\Q}(\lambda^G)$. As $\varepsilon(H,K)(\beta_1u_1)^{-i'}\alpha^{-1}\widehat{E}\alpha(\beta_1u_1)^{j'}\in \operatorname{Cen}_{\mathcal{A}}(\mathcal{B})$ and $ze_{\Q}(\lambda^G)\in\mathcal{B},$ these two commute. Also, $\beta_iu_i,~\alpha$ and $\widehat{E}$ commute with $\varepsilon(H,K)$. This gives the following:\para\noindent 
Matrix units of $\mathcal{A}$: $$\{z^{-i}(\beta_1u_1)^{-i'}\alpha^{-1}\widehat{E}\varepsilon(H,K)\alpha(\beta_1u_1)^{j'}z^{j}~|~1\leq i,j\leq 5,1\leq i',j'\leq 2\}.$$
\para \noindent Primitive orthogonal idempotents of $\mathcal{A}$: $$\{z^{-i}(\beta_1u_1)^{-j}\alpha^{-1}\widehat{E}\varepsilon(H,K)\alpha(\beta_1u_1)^{j}z^{i}~|~1\leq i\leq 5,1\leq j\leq 2\}.$$~~\qed

\bibliographystyle{amsplain}
\bibliography{FinalBG3}
\end{document}